\documentclass[11pt,leqno]{article}
\usepackage{amssymb,amsfonts}
\usepackage{amsmath,amsthm,amsxtra}
\usepackage{epsfig}
\usepackage{slashed}
\usepackage{color}
\usepackage{verbatim}

\newcommand{\ch}{{\rm ch \,}}

\newcommand{\End}{\mathop{\rm End }}

\newcommand{\tr}{\mathrm{tr} \, }

\newcommand{\vac}{|0\rangle}

\renewcommand{\O}{\mathcal{O}}

\newcommand{\CC}{\mathbb{C}}

\newcommand{\ZZ}{\mathbb{Z}}

\newcommand{\fb}{\mathfrak{b}}

\newcommand{\fg}{\mathfrak{g}}
\newcommand{\fh}{\mathfrak{h}}

\newcommand{\fn}{\mathfrak{n}}

\newcommand\bl{(\, . \, | \, . \, )}

\renewcommand{\tilde}{\widetilde}
\renewcommand{\hat}{\widehat}

\makeatletter
\renewcommand\section{\@startsection {section}{1}{\z@}%
                                   {-3.5ex \@plus -1ex \@minus -.2ex}%
                                   {2.3ex \@plus.2ex}%
                                   {\normalfont\large\bfseries}}
\renewcommand\subsection{\@startsection{subsection}{2}{\z@}%
                                     {-3.25ex\@plus -1ex \@minus -.2ex}%
                                     {0ex \@plus .0ex}%
                                     {\normalfont\normalsize\bfseries}}

\@addtoreset{equation}{section}
\makeatother

\newtheorem{theorem}{Theorem}[section]

\newtheorem{lemma}[theorem]{Lemma}

\newtheorem{conjecture}[theorem]{Conjecture}
\newtheorem*{lemma*}{Lemma}

\theoremstyle{remark}
\newtheorem{remark}[theorem]{Remark}
\newtheorem{example}[theorem]{Example}

\makeatletter
\def\@maketitle{\newpage
 \null
 \vskip 2em
 \begin{center}%
 \vskip 3em
  {\Large\bf \@title \par}%
  \vskip 1.5em
  {\normalsize
   \lineskip .5em
   \begin{tabular}[t]{c}\@author
   \end{tabular}\par}%
  \vskip 2em

 \end{center}%
 \par
 \vskip 2.5em}
\makeatother

\renewcommand{\epsilon}{\varepsilon}

\definecolor{light}{gray}{.9}

\newcommand{\half}{\frac{1}{2}}
\newcommand{\thalf}{\tfrac{1}{2}}

\newcommand{\la}{\lambda}
\newcommand{\La}{\Lambda}
\newcommand{\al}{\alpha}

\newcommand{\wg}{\widehat{\fg}}
\newcommand{\wh}{\widehat{\fh}}
\newcommand{\wrh}{\widehat{\rho}}

\newcommand{\LLa}{L(\Lambda)}

\begin{document}

\title{On characers of irreducible  highest weight modules of negative integer level over
  affine Lie algebras}


\author{Victor G.~Kac\thanks{Department of Mathematics, M.I.T, 
Cambridge, MA 02139, USA. Email:  kac@math.mit.edu~~~~Supported in part by an NSF grant.} \ 
and Minoru Wakimoto\thanks{Email: ~~wakimoto@r6.dion.ne.jp~~~~.
Supported in part by Department of Mathematics, M.I.T.}}

\maketitle
\begin{center}
\emph{To the memory of Bertram Kostant}
\end{center}

\setcounter{section}{-1}

\section*{Abstract} We prove a character formula for the irreducible modules from the category $\O$ over the simple affine vertex algebra of type $A_n$ and $C_n$ ($n\geq 2$) of level $k=-1$. We also give
a conjectured character formula for types $D_4$, $E_6$, $E_7$, $E_8$ and levels
$k=-1,...,-b$, where $b=2,3,4,6$ respectively.

\section{Introduction}
Let $ \fg $ be a simple finite-dimensional Lie algebra over $ \CC, $ and let $ \bl $ be the invariant symmetric bilinear form on $ \fg,  $ normalized by the condition $ (\al|\al) = 2 $ for a long root $ \al.  $ Recall that the affine Lie algebra $ \wg, $ associated to $ \fg $ is the infinite-dimensional Lie algebra over $ \CC $
\begin{equation}\label{0.01}
\wg = \fg [t, t^{-1}] \oplus \CC K \oplus \CC d
\end{equation}
with the following commutation relations $ (a, b \in \fg, \ m,n \in \ZZ) $:
\begin{equation}\label{0.02}
[at^m, bt^n] = [a,b]t^{m+n} + m \delta_{m,-n} (a|b)K, \quad [d, at^m] = m\, at^m, \quad [K, \wg] = 0.
\end{equation}
The form $ \bl $ extends from $ \fg $ to a non-degenerate invariant symmetric bilinear form on $ \wg $ by
\begin{equation}\label{0.03}
(at^m|bt^n) = \delta_{m,-n} (a|b), \quad (\fg[t, t^{-1}]| \CC K + \CC d) = 0, \quad (K|K) = (d|d) = 0, \quad (K|d) = 1. 
\end{equation}

Choosing a Cartan subalgebra $ \fh $ of $ \fg, $ one defines the associated Cartan subalgebra of $ \wg: $
\begin{equation}\label{0.04}
\wh = \fh + \CC d +\CC K. 
\end{equation}
The restriction of the bilinear form $ \bl $ from $ \wg $ to $ \wh $ is non-degenerate, and we identify $ \wh $ with $ \wh^* $ and $\fh$ with $\fh^*$ using this form. Then $ d $ and $ K $ are identified with elements, traditionally denoted by $ \La_0 $ and $ \delta $ respectively. We denote $q=e^{-\delta}$. 

Choosing a Borel subalgebra $\fb = \fh + \fn_+ $ of $ \fg, $ one defines the corresponding Borel subalgebra of $ \wg: $
\[ \hat{\fb} = \wh + \fn_+ + \bigoplus\limits_{n>0} \fg t^n . \]
Given $ \La \in \wh^* $, one extends it to a linear function on $ \hat{\fb} $ by zero on all other summands. 
Then there exists a unique irreducible $ \fg $-module $ \LLa, $ which admits an eigenvector
$v_\La$ of $ \hat{\fb} $ with weight $ \La. $ Since $ K $ is a central element of $ \wg,  $ it is represented on $ \LLa $ by a scalar $ \La (K), $ called the level of $ \LLa $ (and of $ \La $).

Let $ \al_1, \ldots, \al_\ell $ be simple roots of $ \fg $, $\theta$ be the highest root, and $ \bar{ \La}_1, \ldots, \bar{ \La}_\ell $ be its fundamental weights, i.e. $ (\bar{ \La}_i | \al^\vee_j) = \delta_{ij}, $ where $ \al^\vee = 2 \al / (\al | \al). $ Then $ \al_0 = \delta - \theta, \al_1, \ldots, \al_{\ell}$ are simple roots of $ \wg, $ and the fundamental weights $\La_i$ of $ \wg $ are defined by 
 $ (\La_i | \al^\vee_j) = \delta_{ij}, \quad \La_i (d) = 0, \ i,j = 0, 1, \ldots, \ell. $
  Any $ \La \in \wh^* $ can be uniquely written in the form 
 \begin{equation}\label{0.05}
 \La = \sum_{i=0}^{\ell} m_i \La_i + a \delta, \text{ where } m_i, a \in \CC. 
 \end{equation}
 
 A $ \wg $-module can be ``integrated'' to the corresponding group, hence it is called integrable, iff all $ m_i $ are non-negative integers.  In this case the level of $ \LLa $ is a non-negative integer. 
 
 The character of $ \LLa $ is defined as the following series, corresponding to the weight space decomposition of $ \LLa $ with respect to $ \wh: $
 \[ (\ch \LLa) (h) = \tr_{\LLa} e^h, \ h \in \wh. \]
  This series is convergent in the domain $ \{ h \in \wh \, | \, \al_i (h) > 0, \ i = 0, 1, \ldots, \ell \}. $ Note that, adding $ b \delta $ to $ \La, $ where $ b \in \CC,  $ multiplies the character by $ q^{-b}. $  Thus, $ \ch \LLa $ depends essentially only on the labels $ m_0, \ldots, m_\ell $ of $ \La $ in \eqref{0.05}.
 
 If the $ \wg $-module $L(\La)$ is integrable, its character is given by the Weyl-Kac character formula:
 \begin{equation}\label{0.06}
 \hat{R} \ch L (\La)  = \sum_{w \in \hat{W}} \epsilon (w) w (e^{\La + \wrh})  = \sum_{w \in W} \epsilon (w) w \sum_{\gamma \in Q^\vee} t_\gamma (e^{\La + \wrh}). 
 \end{equation}
  Here $ \hat{R} = e^{\wrh} \prod_{\al \in \hat{\Delta}_+} (1-e^{-\al})^{\hbox{mult}(\al)} $ is the affine Weyl denominator, $ \wrh $ is the affine Weyl vector:
 \begin{equation}\label{0.07}
 \wrh = \rho + h^\vee \La_0,
 \end{equation}
 where $ \rho  $ is the Weyl vector for $ \fg $ and $ h^\vee $ is the dual Coxeter number ($ =\half $ the eigenvalue on $ \fg $ of the Casimir element). Furthermore, $ \hat{W} = W \ltimes \{ t_\al \, | \, \al \in Q^\vee \} $ is the affine Weyl group, where $ W $ is the Weyl group of $ \fg$, $\epsilon(w) = \det_{\hat{\fh}^*} w$, $ Q^\vee=\sum_{i=1}^\ell\ZZ \al_i^\vee  $ is the coroot lattice of $ \fg, $ and the translation $ t_\gamma \in \End \hat{\fh}^* $ for $ \gamma \in \fh^* $ is defined by 
 \begin{equation}\label{0.08}
 t_\gamma (\la) = \la + \la(K) \gamma - ((\la \, | \, \gamma) + \thalf \la (K) (\gamma \, | \, \gamma)) \delta, \,\,\lambda\in \hat{\fh}^*.
 \end{equation}
 The details of the above discussion may be found in the book \cite{K90}.
 
 In the paper \cite{KW88} a similar character formula was proved for admissible $ \LLa, $ defined by the condition that
 for $ \al \in \hat{\Delta}^{\text{re}}_+ $, the set of real roots of $\fg$,  the number $  (\La + \wrh \, | \, \al^\vee) $ must be a positive integer each time when it is an integer. Admissible $ \wg $-modules include the integrable ones, but exclude, for example the $ \wg $-modules $ L(k\La_0), $ where $ k $ is a negative integer. 
 
It is known for arbitrary (non-critical, i.e. of level $ \neq -h^\vee  $) $ \La  $ that 
\begin{equation}\label{0.09}
 \hat{R} \ch \LLa = \sum_{w \in \hat{W}} c(w) w (e^{\La + \wrh}),
\end{equation}
where $ c(w) $ are integers \cite{KK79}, and that $ c(w) $ can be computed via the Kazhdan-Lusztig polynomials for $ \hat{W}$ [KT00].  However the explicit formulas for the integers $ c(w) $ are unknown in general.
 
In Sections 1 and 2 of the present paper we find explicit character formulas for level $-1$ modules $\LLa$ over
$ \hat{s \ell}_n $ and over $ \hat{sp}_n $ with $n\geq 3$, with highest weights
$ \La = -(1+s) \La_0 + s\La_1 $ and $ \La = -(1+s) \La_0 + s\La_{n-1},
  \ s \in \ZZ$, and $ \La = -(1+s) \La_0 + s\La_1, \ s \in \ZZ_{\geq 0} $
and $ \La = -2\La_0 + \La_2$, respectively (see Theorems \ref{Th1.1}, and \ref{Th2.1}, \ref{Th2.2} respectively). In particular, we compute in both cases the character of $ L(-\La_0), $ which are simple affine vertex algebras of level $ -1. $
(As shown in \cite{AP12} and \cite{AP14}, the above modules are all irreducible modules over these vertex algebras in the category $\O$.)
The main ingredients of the proof are the free field realization of these modules, given in \cite{KW01}, the irreducibility theorems from
\cite{AP14}, \cite{AP12},
and the affine denominator identity for that affine Lie superalgebras
$\hat{ s \ell}_{n|1}$ and $\hat{spo}_{n|2}$, given in \cite{KW94}, \cite{G11}

In Section 3 we indicate a proof, under a certain hypothesis, of an explicit character
formula for certain modules $L(\Lambda)$ of negative integer level over affine Lie algebras,
and conjecture that the hypothesis holds for the affine Lie algebras of Deligne series
$\hat{D}_4$, $\hat{E}_6$, $\hat{E}_7$ and $\hat{E}_8$.

Throughout the paper the base field is $\CC$.

Both of us wish to thank ESI, Vienna, where we began discussion of this paper, for hospitality.
The first named author wishes to thank IHES, where the paper was completed, for perfect working conditions.

\section{Proof of the character formulas for
  $\hat{\fg}$, where $\fg=s\ell_n$}
In this section we prove a character formula for certain highest weight mdules $ \LLa $ of level $ -1 $ over the affine Lie algebra $\hat{s\ell}_n . $ The normalized invariant bilinear form $ \bl $ on $ s\ell_n  $ is the trace form. We choose as its Cartan subalgebra
, as usual, the subalgebra of all diagonal traceless matrices. Then the simple roots of $s\ell_n $ are $ \al_1 = \epsilon_1 - \epsilon_2, \ldots, \al_{n-1} = \epsilon_{n-1} - \epsilon_n , $ where $ \epsilon_1, \ldots, \epsilon_n $ is the standard basis of the dual of all diagonal matrices, and its root (= coroot) lattice is $ Q = \sum_{i=1}^{n-1} \ZZ \al_i. $

We will also use the embedding of $s\ell_n $ in the Lie superalgebra $ s \ell_{n|1}. $ The trace form on $s \ell_n $ extends to the supertrace form $ \bl $ on $ s \ell_{n|1} $, and its Cartan subalgebra embeds in the Cartan subalgebra of $ s \ell_{n|1} $ of supertraceless diagonal matrices. Then the simple roots of $ s \ell_{n|1} $ are $ \al_1, \ldots, \al_{n-1}, \al_n=\epsilon_n-\epsilon_{n+1} $, where $\al_n$ is an odd root. Also, $s \ell_n$ and $ s \ell_{n|1} $ have the same Weyl group, and their dual Coxeter numbers are $ n $ and $ n-1 $ respectively.

\begin{theorem}
\label{Th1.1} Let $ n \geq 3 $, and let $ \LLa $ be an irreducible level $-1$ $\hat{s\ell}_n $-module with highest weight
\[ \La = -(1 + s) \La_0 + s \La_{n-1} \ (\text{resp. }\, = -(1+s)\La_0  + s\La_1),\, s\in \ZZ_{\geq0}. \]
Then the character of $ \LLa $ is given by the following formula:
\begin{equation}
\label{1.01}
\hat{R} \ch \LLa = \sum_{w \in W} \epsilon (w) w \sum_{\substack{\gamma \in Q \\ (\gamma |\bar{ \La}_{n-1} (\text{resp. }\bar{ \La}_1)) \geq 0}} t_\gamma (e^{\La + \hat{\rho}}).
\end{equation}
\end{theorem}

The proof of formula \eqref{1.01} uses the free field construction, given in \cite{KW01}, of the affine Lie superalgebra $ \hat{g \ell}_{m|n} $ of level 1 in a Fock space $ F $ in the case of $ m = 0 $
. Note that in that paper we used the supertrace form, which is equal to the negative of the trace form on $ g\ell_n = g \ell_{0|n}. $ Hence we get a $ \hat{g \ell}_n $-module structure on $ F $ of level $ -1. $
Recall some properties of this module, described in \cite{KW01}.

First, we have the charge decomposition into a direct sum of $ \hat{g \ell}_n $-submodules:
\begin{equation}\label{1.02}
F = \underset{s \in \ZZ}{\bigoplus} F_s.
\end{equation}

Second, there is a Virasoro algebra acting on $ F, $ and leaving all subspaces $ F_s $ invariant, for which all fields $ a(z), \ a \in g\ell_n, $ are primary of conformal weight 1, and each $ F_s $ in \eqref{1.02} has a unique, up to a constant factor, non-zero vector $ | s \rangle $ with minimal $L_0$-eigenvalue (see Section 2 for more details). Moreover this vector is invariant with respect to the Cartan subalgebra of $\hat{ g\ell}_n $ and has the following weight:
\begin{equation}\label{1.03}
\text{weight } |s\rangle = \begin{cases}
-\La_0 - \tfrac{s}{2} \delta + s \epsilon_1 & \text{if } s \in \ZZ_{\geq 0}, \\
-\La_0 + \tfrac{s}{2} \delta + s \epsilon_n & \text{if } s \in \ZZ_{\leq 0}. \\
\end{cases}
\end{equation}

Third, by formula (3.15) from \cite{KW01}, the character of the $\hat{ g \ell}_n $-module $ F $ is given by 
\[ \ch F : = \sum_{s \in \ZZ} x^s \ch F_s = e^{-\La_0} \prod_{j =1}^{n} \prod_{k =1 }^{\infty} (1 - xe^{\epsilon_j} q^{k-\half})^{-1} (1 - x^{-1}e^{-\epsilon_j} q^{k-\half})^{-1}. \]
Letting in this formula $ x = e^{-\epsilon_{n+1}}q^\half, $
we obtain:
\begin{equation}\label{1.04}
e^{-\La_0} \prod_{j =1}^{n} \prod_{k =1 }^{\infty} (1 - e^{\epsilon_j - \epsilon_{n+1}} q^k)^{-1} (1 - e^{-(\epsilon_j -\epsilon_{n+1} )} q^{k-1})^{-1} = \sum_{s \in \ZZ} e^{-s \epsilon_{n+1}} q^{\frac{s}{2}} \ch F_s.
\end{equation}

It was proved in \cite{AP14} that all $ \hat{g \ell}_n $-modules $ F_s $ are irreducible, provided that $ n \geq 3. $
Therefore, using \eqref{1.03}, we see that $ F_s = V (\la^{(s)}) \otimes L (\La^{(s)}) $, where $ V (\la^{(s)}) $ is an irreducible $ \hat{g\ell}_1 $-module with highest weight $ \la^{(s)} \in \CC \sum_{i = 1}^{n} \epsilon_i + \CC \delta $ and $ L(\La^{(s)}) $ is an irreducible $ \hat{s \ell}_n $-module with highest weight 
\[ \La^{(s)} = -\La_0 + s \bar{ \La}_1 \ (\text{resp. } -\La_0 - s \bar{ \La}_{n-1}) \in \fh^* \text{ if } s \geq 0 \ (\text{resp. } s \leq 0), \]
where $ \la^{(s)} \oplus \La^{(s)} = $ weight $ |s \rangle $. 

Hence, using that $\La_i=\La_0 +\bar{\La}_i$, we obtain that the character of $ F_s $
is given by
\begin{equation}\label{1.05}
\varphi (q) \ch F_s = \begin{cases}
q^{\frac{s}{2}} e^{s(\epsilon_1 - \bar{\La}_1)} \ch L (-(1+s)\La_0 + s\La_1) & \text{if } s \in \ZZ_{\geq 0}, \\
q^{-\frac{s}{2}} e^{s(\epsilon_n + \bar{\La}_{n-1})} \ch  L (-(1-s)\La_0 - s\La_{n-1}) & \text{if } s \in \ZZ_{\leq 0}. \\
\end{cases}
\end{equation}
Here and further $\varphi(q)=\prod_{n=1}^\infty(1-q^n)$.
Substituting \eqref{1.05} in the RHS of \eqref{1.04}, we obtain:
\begin{equation}
\label{1.06}
\begin{aligned}
\text{LHS of \eqref{1.04}} & = \frac{1}{\varphi (q)} ( \sum_{s>0} e^{s(\epsilon_1 - \epsilon_{n+1})} q^s e^{-s \bar{\La}_1} \ch L (-(1+s)\La_0 + s\La_1) . \\ 
&  + \sum_{s\leq 0} e^{s(\epsilon_n - \epsilon_{n+1})} e^{s\bar{\La}_{n-1}} \ch L (-(1-s)\La_0 - s\La_{n-1})  ). \\
\end{aligned}
\end{equation}
Next, we embed the Lie algebra $ s \ell_n $
in the Lie superalgebra $ s\ell_{n|1} $ as described above.
We extend this embedding to the affine Lie (super)algebras $ \hat{s \ell}_n \hookrightarrow \hat{s \ell}_{n|1}. $ Then identitiy \eqref{1.06} can be rewritten as follows:
\begin{equation}\label{1.07}
\begin{aligned}
& e^{-\La_0} \varphi (q) \prod_{j =1}^n \prod_{k =1 }^\infty (1 - e^{\al_j + \cdots + \al_n} q^k )^{-1} ( 1 - e^{-(\al_j + \cdots + \al_n)} q^k)^{-1} \\
& = \sum_{s>0} e^{s(\al_1 + \cdots + \al_n)} q^s e^{-s\bar{\La}_1} \ch L (-(1+s)\La_0 + s \La_1) + \sum_{s\leq 0} e^{s \al_n} e^{s \bar{\La}_{n-1}} \ch L (-(1-s)\La_0 - s\La_{n-1}). \\
\end{aligned}
\end{equation}

We denote the $ 0 $-th fundamental weight and the Weyl vector for $ \hat{s \ell}_{n|1} $ by $ \La'_0 $ and $ \wrh' $ respectively. Then, by \eqref{0.07}, we have, when restricted to $ \hat{s \ell}_n: $
\begin{equation}\label{1.08}
\La'_0 = \La_0 \text{ and } \wrh' = \wrh - \La_0.
\end{equation}
Recall the formulas for the Weyl denominator $ \hat{R} $ for $ \hat{s \ell}_n $ and the Weyl superdenominator $ \hat{R}' $ for $ \hat{s \ell}_{n|1}: $
\begin{equation}\label{1.09}
  \hat{R} = e^{\wrh} \varphi (q)^{n-1} \prod_{\al \in
    \hat{\Delta}^\text{re}_+}
  (1-e^{-\al}),
\end{equation}
\begin{equation}\label{1.10}
\hat{R}' = e^{\wrh'} \varphi (q)^n \displaystyle\prod_{\al \in \hat{\Delta}^\text{re}_+} (1-e^{-\al})\displaystyle\prod_{j =1}^n \displaystyle\prod_{k =1 }^\infty (1-e^{\al_j + \cdots + \al_n} q^k)^{-1} (1-e^{-(\al_j + \cdots + \al_n)} q^{k-1})^{-1}.
\end{equation}
Multiplying both sides of \eqref{1.07} by $ \hat{R}, $ we obtain, using \eqref{1.09} and \eqref{1.10}:
\begin{equation}\label{1.11}
\begin{aligned}
\hat{R}' &= \sum_{s>0} e^{s(\al_1 + \cdots + \al_n)} q^s e^{-s \bar{\La}_1} \hat{R} \ch L (-(1+s)\La_0 + s\La_1) \\
& + \sum_{s\leq 0} e^{s \al_n} e^{s \bar{\La}_{n-1}} \hat{R} \ch L (-(1-s)\La_0 - s\La_{n-1}).\\
\end{aligned}
\end{equation}

On the other hand, $ \hat{R}' $ can be computed by the superdenominator identity \cite{KW94}, \cite{G11}:
\begin{equation}\label{1.12}
\hat{R}' = \epsilon (w) w \sum_{\gamma \in Q} t_\gamma \frac{e^{\wrh'}}{1 - e^{-\al_n}}.
\end{equation}
Expanding $ t_\gamma \frac{e^{\wrh'}}{1 - e^{-\al_n}} $ in the geometric series in the domain $ |e^{-\al_n}| < 1, \ |q| <1, $ we obtain for $ \gamma \in Q$
(using \eqref{0.08}): 
\[ t_\gamma \frac{e^{\wrh'}}{1 - e^{-\al_n}} = e^{(n-1)\La_0} ( \sum_{\substack{p \geq 0 \\ (\gamma | \al_n)\leq 0 }}- \sum_{\substack{p < 0 \\ (\gamma | \al_n) > 0 }}) e^{\rho + (n-1) \gamma - p \al_n} q^{\frac{n-1}{2} (\gamma| \gamma) + (\rho | \gamma) - p (\al_n|\gamma)}. \]
Using that $ n \al_n = \sum_{i = 1}^{n} \epsilon_i -n\bar{ \La}_{n-1},$ and that $(\al_n|\gamma)=-(\bar{\La}_{n-1}|\gamma)$ for $\gamma\in Q$, we can rewrite this formulas as
\[  t_\gamma \frac{e^{\wrh'}}{1 - e^{-\al_n}} = e^{(n-1) \La_0} ( \sum_{\substack{p \geq 0 \\ (\gamma | \bar{\Lambda}_{n-1})\geq 0 }}- \sum_{\substack{p < 0 \\ (\gamma | \bar{\Lambda}_{n-1}) < 0 }}) e^{-\frac{p}{n} \sum_{i = 1}^{n} \epsilon_i} e^{\rho + (n-1) \gamma + n \bar{ \La}_{n-1}} q^{\frac{n-1}{2} (\gamma|\gamma) + (\rho|\gamma)+ p (\bar{\La}_{n-1} | \gamma)}. \]
Plugging this in \eqref{1.12} and using that $ \sum_{i=1}^{n} \epsilon_i $ is $ W $-invariant, we obtain:
\begin{equation}\label{1.13}
\begin{aligned}
  \hat{R}'  =& e^{(n-1)\La_0} \sum_{\gamma \in Q} ( \sum_{\substack{p \geq 0 \\ (\gamma | \bar{\La}_{n-1})\geq 0 }}- \sum_{\substack{p < 0 \\ (\gamma | \bar{\La}_{n-1}) < 0 }}) e^{-p(\al_n + \bar{ \La}_{n-1})} q^{\frac{n-1}{2} (\gamma|\gamma) + (\rho|\gamma) +
    p (\bar{ \La}_{n-1} | \gamma)} \\
&\times \sum_{w \in W} \epsilon (w) w e^{\rho +(n-1) \gamma + p \bar{ \La}_{n-1}}. \\
\end{aligned}
\end{equation}
Thus we obtain the identity
\[ \text{RHS of \eqref{1.11}} = \text{ RHS of \eqref{1.13}}. \]
Comparing the coefficient of $ e^{-p\al_n} $ for $ p \geq 0 $ in this identity, we obtain:
\[ 
\hat{R} \ch L (-(1+p)\La_0 + p \La_{n-1}) = e^{(n-1) \La_0}  \sum_{\substack{\gamma \in Q \\ (\gamma | \bar{\La}_{n-1}) \geq 0}} q^{\frac{n-1}{2} (\gamma | \gamma)+(\rho | \gamma)+ p (\bar{ \La}_{n-1} | \gamma)}\] 
\[\times \sum_{w \in W} \epsilon (w) w e^{\rho + (n-1) \gamma + p \bar{ \La}_{n-1}} 
=  \sum_{w \in W} \epsilon (w) w \sum_{\substack{\gamma \in Q \\ (\gamma | \bar{\La}_{n-1}) \geq 0}} t_\gamma e^{-(1+p)\La_0 + p \La_{n-1} + \wrh}.\] 
This establishes formula \eqref{1.01} for $ \La = - (1+s) \La_0 + s \La_{n-1}$.
Formula for $ \La = -(1+s) \La_0 + s \La_1 $ follows by the involution of the Dynkin diagram of $ \hat{s \ell}_n $ which keeps the $ 0 $th node fixed. 
\begin{remark}
\label{rem1.2}
Let $ \La = -(1+s) \La_0 + s \La_{1} $ and let $ \LLa = \bigoplus\limits_{j \in \ZZ_{\geq 0}} \LLa_j $ be the eigenspace decomposition of $ \LLa $ with respect to $ -d. $ Let $ \dim_q \LLa = \sum_{j \geq 0} (\dim \LLa_j)q^j $ be the ``homogeneous'' $ q $-dimension of $ \LLa $. Dividing both sides of the last equality in the proof of \eqref{1.01} by $ e^{(n-1)\La_0} R, $ where $ R $ is the Weyl denominator of $ \fg,  $ and letting all elements of $ \fh^* $ equal 0, we obtain by the usual argument:
\[ \varphi(q)^{\dim \fg} \dim_q \LLa = \sum_{\substack{\gamma \in Q \\ (\bar{\La}_{1}| \gamma) \geq 0}} \dim (s \bar{\La}_{1} + (n-1) \gamma) q^{\frac{n-1}{2} (\gamma| \gamma) + (s\bar{\La}_{1} + \rho | \gamma)}, \]
where $ \dim (\la)=\prod_{\al \in \Delta_+}(\la+\rho|\al)/(\rho|\al) $ is the expression of the Weyl dimension formula for $ \fg. $
\end{remark}
\begin{remark}
\label{rem1.3}
For $ \hat{s\ell}_2 $ the characters of the above modules are
very easy (see, e.g. \cite{KW88}):
$ e^{\La_0} \hat{R} \ch L (-(1+s)\La_0 + s\La_1) = 1 - e^{-(s+1)\al_1}, s\in\ZZ_{\geq 0}.$
\end{remark}

\section{Proof of the character formulas for $\hat{\fg}$,
 where $\fg=sp_n$}

In order to prove the character formulas for $ \hat{sp}_n, $ where $ n \geq 4,  $ even, we use the embedding of $ sp_n $ in $ s \ell_n, $ obtained as the fixed point subalgebra of the involution $ \sigma, $ corresponding to the flip of the Dynkin diagram of $ s \ell_n. $ The normalized invariant bilinear form $ \bl $ on $ sp_n $ is again the trace form. Let $ n' = \frac{n}{2}. $

We will also use the embedding of the Lie algebra $ sp_n $ in the Lie superalgebra $ spo_{n|2}. $ The simple roots of the latter are $ \al_*, \al_1, \ldots, \al_{n'-1}, \al_{n'},$ where
$\al_*$ is an odd root, and $ \al_1, \ldots, \al_{n'} $ are the simple roots of the former. We let $ Q^\vee = \sum_{i=1}^{n'} \ZZ \al^\vee_i $ be the coroot ($ \neq $ root) lattice of $ sp_n. $

In this section we prove the following two theorems. 

\begin{theorem}
\label{Th2.1}
Let $ n \geq 4 $ be even, and let $ \LLa $ be an irreducible $ \hat{sp}_n $-module with highest weight $ \La. $
\begin{enumerate}
\item[(a)] If $ \La = -(1+s) \La_0 + s\La_1,  $ where $ s \in \ZZ_{\geq 1}, $ then
\[ \ch \LLa = \hat{R}^{-1} \sum_{w \in W} \epsilon (w) w \sum_{\substack{\gamma \in Q^\vee \\ (\gamma |\bar{\La}_{1} ) \geq 0}}  t_\gamma (e^{\La + \wrh})\]
\item[(b)] If $ \La = -\La_0, $ then 
  \[ \ch \LLa = \half (  \hat{R}^{-1} \sum_{w \in W} \epsilon (w) w \sum_{\substack{\gamma \in Q^\vee \\ (\gamma |\bar{\La}_{1} ) \geq 0}}  t_\gamma (e^{\La + \wrh}) + e^{-\La_0}
  \frac{\varphi (q^2)}{\varphi(q)} \prod_{\al \in \Delta_\ell} \prod_{k \in \ZZ_{\text{odd}>0}} (1-e^\al q^k)^{-1}  ),  \]
    where $\Delta_\ell$ is the set of long roots of $sp_n$.
\item[(c)] If $ \La = -2 \La_0 + \La_2, $ then $ \ch \LLa $ is obtained by an expression, obtained from (b) by dividing by $ q $ and replacing plus by minus between the summands.
\end{enumerate}
\end{theorem}

\begin{theorem}
\label{Th2.2}
The characters of the $ \hat{sp}_n $-modules $ L(-\La_0) $ and $ L(-2\La_0 + \La_2) $ can be written in the form \eqref{0.09} as follows:
\begin{enumerate}
\item[(a)]
\[ \hat{R} \ch L(-\La_0) = \sum_{w \in W} \epsilon (w) w \sum_{\substack{\gamma \in Q^\vee \\ (\gamma |\bar{\La}_{1} ) \geq 0 \\ (\gamma | \bar{ \La}_{n'}) \in 2 \ZZ}}  t_\gamma (e^{\La + \wrh}). \] 
\item[(b)] 
\[ \hat{R} \ch L(-2\La_0 + \La_2) = \sum_{w \in W} \epsilon (w) w \sum_{\substack{\gamma \in Q^\vee \\ (\gamma |\bar{\La}_{2}-\bar{ \La}_1 ) \geq 0 \\ (\gamma | \bar{ \La}_{n'}) \in 2 \ZZ}}  t_\gamma (e^{\La + \wrh}). \]
\end{enumerate}
\end{theorem}

First letting in the character $ \ch F $ of the $ \hat{g \ell}_n $-module $ F, $ considered in Section 1, $  x = e^{\al_* - \epsilon_{n'}} q^\half, $ and restricting the module $ F $ to $ \hat{sp}_n, $ we obtain, cf. \eqref{0.04}:
\begin{equation}\label{2.01}
\begin{aligned}
&e^{-\La_0} \prod_{k=1}^{\infty} (1-e^{\al_*}q^k) (1-e^{-\al_*}q^{k-1})  \prod_{k =1 }^{\infty} \prod_{i =1}^{n'} (1-e^{\al_* +\al_1 + \cdots + \al_i}q^k) (1-e^{-\al_* -\al_1 - \cdots -\al_i}q^{k-1})  \\
&\times \prod_{k =1 }^{\infty} \prod_{i =1}^{n'-1} (1-e^{\al_* +\al_1 + \cdots + 2\al_i + \cdots + 2 \al_{n'-1} + \al_{n'}}q^k) (1-e^{-\al_* -\al_1 - \cdots - 2\al_i - \cdots - 2 \al_{n'-1} -\al_{n'}}q^{k-1}) \\
&= \sum_{s \in \ZZ} e^{s(\al_* - \epsilon_{n'})} q^{\frac{s}{2}} \ch \, F_s \big{|}_{\hat{sp}_n}.
\end{aligned}
\end{equation}
Next, we embed the Lie algebra $ sp_n $ in the Lie superalgebra $ spo_{n|2} $ as described above. We denote the $ 0 $-th fundamental weight and the Weyl vector for $ \hat{spo}_{n|2} $ by $ \La_0' $ and $ \wrh' $ respectively. Then, by \eqref{0.07}, we have again \eqref{1.08}, when restricted to $ \hat{sp}_n. $

Denote by $ \hat{R} $ and $ \hat{R}' $ the Weyl denominator and Weyl superdenominator for $ \hat{sp}_n $ and $ \hat{spo}_{n|2} $ respectively. Then we have by \eqref{1.08} for $ sp_n: $ 
\begin{equation}\label{2.02}
\hat{R}' = \varphi (q) \hat{R} \times (\text{LHS of \eqref{2.01}}).
\end{equation}
On the other hand, the superdenominator identity for $ \hat{spo}_{n|2} $ reads \cite{KW94}, \cite{G11}:
\begin{equation}\label{2.03}
\hat{R}' = \sum_{w \in \hat{W}} \epsilon (w) w \frac{e^{\wrh'}}{1-e^{-\al_*}},
\end{equation}
where $ \hat{W} = W \ltimes \{  t_\gamma \ | \ \gamma \in Q^\vee    \} $ is the Weyl group of $ \hat{sp}_n. $

Using \eqref{2.01}--\eqref{2.03} and applying the same argument as in Section 1, we obtain for each $ s \in \ZZ_{\geq 0}: $
\begin{equation}\label{2.04}
\varphi(q)\hat{R} \,  (\ch F_s) \big{|}_{\hat{sp}_n} = \sum_{w \in W} \epsilon (w) w \sum_{\substack{\gamma \in Q^\vee \\ (\gamma |\bar{\La}_{1} ) \geq 0}} t_\gamma (e^{\La + \wrh}).
\end{equation}
This proves claim (a) of Theorem \ref{Th2.1}, since, due to \cite{AP12}, the restriction of the $\hat{s\ell}_n$-module $L(-\La_0+s\bar{\La}_1)$ to $\hat{sp}_n$ is irreducible for $s>0$.

In order to prove claims (b) and (c), we need to study the module $ F $ more carefully. Recall that $ F $ is the unique irreducible module over the Clifford algebra $ Cl $ with generators $ \varphi^{(i)}_k $ and $ \varphi^{(i)*}_k,  \ i = 1, \ldots, n, \ k \in \half + \ZZ, $ with relations $ \varphi^{(i)*}_{-k} \varphi^{(i)}_{k} - \varphi^{(i)}_{k} \varphi^{(i)*}_{-k} = 1 $ and $ =0 $ in the rest of the cases, which admits a non-zero vector $ \vac, $ such that $ \varphi^{(i)}_k \vac = 0 =\varphi^{(i)*}_k \vac  $ for all $ k>0, \ i = 1, \ldots, n. $ The charge decomposition \eqref{1.02} is defined by letting charge $ \vac = 0, $ charge $ \varphi^{(i)}_k = 1 = -\text{charge} \, \varphi^{(i)*}_k  $.

Note that the algebra $ Cl $ carries an involution $ \sigma,  $ defined by 
\[ \sigma (\varphi^{(i)}_k) = (-1)^i \varphi^{(n+1-i)*}_k \]
This involution induces an involution of the space $ F_0,  $ denoted again by $ \sigma, $ letting $ \sigma \vac = \vac, $ so that we have its eigenspace decompositions
\begin{equation}\label{2.05}
F_0 = F^1_0 \oplus F^{-1}_0.
\end{equation}
The space $ F_0  $ is spanned by monomials
\begin{equation}\label{2.06}
v = \varphi^{(i_1)}_{-k_1} \cdots \varphi^{(i_m)}_{-k_m}\varphi^{(j_1)*}_{-l_1} \cdots \varphi^{(j_m)*}_{-l_m} \vac.
\end{equation}
Since 
\[ \sigma (v) =  (-1)^{ \sum_{p=1}^{m} (i_p + j_p)+m} \varphi^{(n+1-j_1)}_{-l_1} \cdots \varphi^{(n+1-j_m)}_{-l_m} \varphi^{(n+1-i_1)*}_{-k_1} \cdots \varphi^{(n+1-i_m)*}_{-k_m} \vac,  \]
we see that, if $ v \in F^{-1}_0, $  we have $ (1\leq p \leq m): $
\begin{equation}\label{2.07}
n+1 - j_p = i_p, \quad n+1 - i_p = j_p, \quad l_p = k_p;
\end{equation}
\begin{equation}\label{2.08}
\sum_{p=1}^{m} (i_p + j_p) + m \equiv 1 \mod 2.
\end{equation}
By \eqref{2.07} we have 
$\sum_{p=1}^{m} i_p = \sum_{p=1}^{m} (n+1) - \sum_{p=1}^{m} j_p.$
Therefore
$\sum_{p=1}^{m} (i_p + j_p) = m \mod 2$,
which contradicts \eqref{2.08}. Hence $ F^{-1}_0 $ contains no monomials \eqref{2.06}.

Thus, for a monomial \eqref{2.06} we have: either $ \sigma (v) = v,  $ or $ v $ and $ \sigma (v) $ are linearly independent. Denote by $ F^\sharp_0 $ the subspace of $ F_0 $ spanned by monomials fixed by $ \sigma. $ From the above discussion we obtain:
\begin{equation}\label{2.09}
\ch F^{\pm 1}_0 = \half (\ch \, F_0 \pm \ch \, F^\sharp_0).
\end{equation}

Recall the construction of the representation of $ \hat{ g\ell}_n $ of level $ -1 $ in $ F $ \cite{KW01}. Let $ \varphi^{(i)} (z) = \sum_{k \in \half + \ZZ} \varphi^{(i)}_k z^{-k-\half}, \ \varphi^{(i)*} (z) = \sum_{k \in \half + \ZZ} \varphi^{(i)*}_k z^{-k-\half}. $ Then
\begin{equation}\label{2.10}
e_{ij} (z) \mapsto : \varphi^{(i)} (z) \varphi^{(j)*} (z) :, \ K \mapsto -1, \ d \mapsto -L_0 
\end{equation}
defines a representation of $ \hat{ g\ell}_n $ in $ F$ of level $ -1 $ (preserving \eqref{1.02}). Here
\[ \sum_{k \in \ZZ} L_k z^{-k-2} = \half \quad \sum_{j=1}^{n} (: \partial \varphi^{(j)} (z) \varphi^{(j)*} (z): - : \partial \varphi^{(j)*} (z) \varphi^{(j)} (z) :) \]
is the representation in $ F $ of the Virasoro algebra. In particular, the Heisenberg subalgebra $ H $ of $ \hat{ g\ell}_n $ is represented in $ F $ as 
\[\sum_{k \in \ZZ} (I_n t^k)z^{-k-1} \mapsto \sum_{k \in \ZZ} H_k z^{-k-1} = \sum_{i=1}^{n} : \varphi^{(i)} (z) \varphi^{(i)*} (z):.\] 
Since the $\hat{gl}_n$-module $F_0$ is irreducible \cite{AP14}, we have the following decomposition of it as an $H\oplus \hat{sl}_n$-module
\begin{equation}\label{2.11}
F_0=V\otimes L(-\Lambda_0),
\end{equation}
where $ V $ is an irreducible $ H $-module with highest weight vector $ \vac, $ i.e. $ (I_nt^k) \vac = 0 $ for $ k \geq 0. $  
We obviously have:
\[ V = \CC \left[ H_{-k} \mid k \in \ZZ_{>0}\right] \vac \text{ and } \sigma (H_k)  = -H_k.\]
Hence, in particular $ V $ is $ \sigma $-invariant, so that we have the eigenspace decomposition with respect to $ \sigma: V = V^1 \oplus V^{-1}. $
Note that the action of $ \hat{sp}_n \subset \hat{ g\ell}_n $ on $ F_0 $ commutes with the action of $ \sigma $ on $ F_0,  $ hence  both $ F^1_0 $ and $ F^{-1}_0 $ are $ \hat{sp}_n $-modules. Moreover, due to \cite{AP12}, the $ \hat{s \ell}_n $-module $ L(-\La_0), $ restricted to $ \hat{sp}_n, $ is a direct sum of two irreducible modules, with highest weights $ -\La_0 $ and $ -\La_0 + \bar{\La}_2 \mod \CC \delta. $ But it is easy to see that 
\[ ( \varphi^{(1)}_{-\half} \varphi^{(n-1)*}_{-\half} + \varphi^{(2)}_{-\half} \varphi^{(n)*}_{-\half}) \vac \]
is a singular vector for $ \hat{sp}_n, $ and its weight is $ -\La_0 + \bar{ \La}_2 - \delta. $
Thus we obtain
\begin{lemma}
\label{Lem2.1}
As an $ H \oplus \hat{sp}_n $-module, one has 
\[ F^1_0 \simeq V^1 \otimes L(-\La_0)+V^{-1}\otimes L(-\Lambda_0+\Lambda_2 -\delta)\,, \quad F^{-1}_0 = V^{-1} \otimes L(-\Lambda_0)+V^1\otimes L(-\La_0 + \bar{ \La}_2 - \delta). \]
\end{lemma}
It is easy to see that
\begin{equation}\label{2.12}
\ch V^{\pm 1} = \half \left( \frac{1}{\varphi (q)} \pm \frac{\varphi (q)}{\varphi (q^2)} \right),
\end{equation}
hence we have
\begin{equation}\label{2.13}
(\ch V^1)^2 - (\ch V^{-1})^2 = \frac{1}{\varphi (q^2)}.
\end{equation}

Next, we obviously have:
\[ F^\sharp_0 = \CC \left[ \varphi^{(i)}_{-k}  \varphi^{(n+1-i)*}_{-k} \, \middle| \, 1 \leq i \leq n, \ k \in \half + \ZZ_{\geq 0} \right] \vac,  \]
hence
\[ \ch F^\sharp_0 = e^{-\La_0} \prod_{k \in \ZZ_{\text{odd}>0}} \prod_{i=1}^{n} (1-e^{\epsilon_i-\epsilon_{n+1-i}} q^k)^{-1}. \]
It follows that
\begin{equation}\label{2.14}
\ch F^\sharp_0 \big{|}_{\hat{sp}_n} = e^{-\La_0} \prod_{k \in \ZZ_{\text{odd}>0}} \prod_{\al \in \Delta_\ell} (1-e^\al q^k)^{-1}.
\end{equation}

Now we are able to complete the proofs of claims (b) and (c) of Theorem \ref{Th2.1}. By Lemma \ref{Lem2.1} we have:
\begin{equation}\label{2.15}
\ch V^{\pm1} \ch L(\La_0) + \ch V^{\mp1} \ch (-\La_0 + \bar{ \La}_2 -\delta) = \ch F^{\pm1}_0.
\end{equation}
From \eqref{2.09}, \eqref{2.13} and \eqref{2.15} we obtain:
\[ \frac{1}{\varphi (q^2)} \ch L(\La_0) = \half \left( \ch V^1 - \ch V^{-1}  \right) \ch F_0 \big{|}_{\hat{sp}_n} + \half (\ch V^1 + \ch V^{-1}) \ch F^\sharp_0 |_{\hat{sp}_n}.\]
Now claim (b) follows from \eqref{2.14}. Claim (c) follows from Lemma \ref{Lem2.1} and claims (a), (b).

Next, we turn to the proof of Theorem \ref{Th2.2}. First from the denominator identity of $ A^{(2)}_{2n'-1} $ we deduce the following lemma. 

\begin{lemma}
\label{Lem2.2}
Let $ M = \{  \gamma \in Q^\vee \, | \, (\gamma | \bar{ \La}_n) \in 2 \ZZ \} $. Then
\[ e^{-\La_0} \frac{\varphi (q^2)}{\varphi (q)} \hat{R} \prod_{\substack{\al \in \Delta_\ell \\ k \in \ZZ_{\text{odd}>0}}} (1-e^\al q^k)^{-1} = \sum_{w \in W} \epsilon (w) w \sum_{\gamma \in M} t_\gamma (e^{n' \La_0 + \rho}). \]
\end{lemma}
Using this lemma, we can rewrite the character formulas, given by Theorem 2.1(b) and (c) as follows
\begin{equation}\label{2.16}
\hat{R} \, \ch L(-\La_0) = \half \sum_{w \in W} \epsilon (w) w ( \sum_{\substack{\gamma \in Q^\vee \\ (\gamma|\bar{ \La}_1) \geq 0}} + \sum_{\substack{\gamma \in Q^\vee \\ (\gamma|\bar{ \La}_{n'}) \in 2 \ZZ}})\, t_\gamma (e^{n'\La_0 + \rho}).
\end{equation}
\begin{equation}\label{2.17}
\hat{R} \, \ch (-\La_0 + \bar{ \La}_2) = \frac{1}{2q} \sum_{w \in W} ( \sum_{\substack{\gamma \in Q^\vee \\ (\gamma|\bar{ \La}_1) \geq 0}} - \sum_{\substack{\gamma \in Q^\vee \\ (\gamma|\bar{ \La}_{n'}) \in 2 \ZZ}} )\, t_\gamma (e^{n'\La_0 + \rho}).
\end{equation}
In order to rewrite these formulas into a nicer form we introduce a different $ \ZZ $-basis of $ Q^\vee $:
\[ \gamma_i = \al^\vee_i + \cdots + \al^\vee_{n'}, \ i = 1, \ldots, n'. \]
Then, letting $ \gamma = \sum_k j_k \gamma_k, $ we have:
\begin{equation}\label{2.18}
(\gamma | \bar{ \La}_1) = j_1, \quad (\gamma|\bar{ \La}_2 - \bar{\La}_1) = j_2, \quad (\gamma| \bar{ \La}_{n'}) = \sum_k j_k.
\end{equation}
Using that $ \left( -\La_0 + \wrh \, \middle| \, \delta - \theta \right) = 0, $ we obtain
\begin{lemma}
\label{Lem2.3}
For $ \Omega \subset \ZZ^{n'} $ let
\[ \Omega' = \left\{ (-j_1-1, j_2, \ldots, j_{n'}) \, \middle| \, (j_1, \ldots, j_{n'}) \in \Omega \right\} .\]
Then 
\[ \sum_{w \in W} \epsilon (w) w \sum_{(j_1, \ldots, j_{n'}) \in \Omega} t_{\sum_k j_k \gamma_k} (e^{n \La_0 + \rho}) =- \sum_{w \in W} \epsilon (w) w \sum_{(j_1, \ldots, j_{n'}) \in \Omega'} t_{\sum_k j_k \gamma_k} (e^{n \La_0 + \rho}). \]
\end{lemma}
Introduce the following shorthand notation:
\[ [\text{condition } (*) \text{ on } \gamma] : = \sum_{w \in W} \epsilon (w) w \sum_{\substack{\gamma \in Q^\vee \\  \gamma \text{ satisfies (*)}}} t_\gamma (e^{n'\La_0 + \rho}).   \]
Applying Lemma \ref{Lem2.3} to the set
$ \Omega = \{ (j_1, \ldots , j_n) \in \ZZ^{n'} \, | \, j_1 \geq 0, \sum_k j_k \in \ZZ_{\text{odd}}   \} $, we obtain in this notation:
\begin{equation}\label{2.19}
\left[ (\gamma| \bar{ \La}_1) \geq 0, (\gamma | \bar{ \La}_{n'}) \in 1 + 2 \ZZ \right] =- \left[ (\gamma | \bar{ \La}_1) < 0, (\gamma| \bar{ \La}_{n'}) \in 2 \ZZ  \right].
\end{equation}
In the above notation, formula \eqref{2.16} becomes:
\[ \hat{R} \, \ch L(-\La_0) = \half (  [(\gamma| \bar{ \La}_1) \geq 0] + [(\gamma| \bar{ \La}_{n'}) \in 2 \ZZ ]) . \]
Using \eqref{2.19}, this completes the proof of claim (c) of Theorem \ref{Th2.2}.

Likewise in the above notation formula \eqref{2.17} becomes:
\[ \hat{R} \, \ch L (-\La_0 + \bar{ \La}_2)  = \frac{1}{2q}( [(\gamma|\bar{ \La}_1) \geq 0]- [(\gamma|\bar{ \La}_{n'}) \in 2 \ZZ]). \]
Using \eqref{2.19} this can be rewritten as 
\begin{equation}\label{2.20}
\hat{R} \, \ch L(-\La_0 + \bar{ \La}_2) = -\frac{1}{q} \sum_{w \in W} \epsilon (w) w \sum_{\substack{\gamma \in Q^\vee \\ (\gamma|\bar{ \La}_1) < 0 \\ (\gamma | \bar{ \La}_{n'}) \in 2 \ZZ  }} t_\gamma (e^{n' \La_0 + \rho}).
\end{equation}
In order to rewrite this formula further we need the following properties of roots and weights of $ sp_n, $ which are straightforward. 
\begin{lemma}
\label{Lem 2.4}
The weight $ \la := \bar{ \La}_2 $ satisfies the following properties:
\begin{enumerate}
\item[(a)] $ \la $ is a positive short root, given by
\[ \la = \al_1 + 2 (\al_1 + \cdots + \al_{n'-1}) + \al_{n'} = \half (\gamma_1 + \gamma_2). \]
\item[(b)] $ (\la| \gamma_i) =1 $ if $ i = 1, 2, $ and $ =0 $ otherwise. 
\item[(c)] \[ r_\la ( \sum_{k=1}^{n'} j_k \gamma_k ) = -j_2 \gamma_1 - j_1 \gamma_2 + \sum_{k=3}^{n'} j_k \gamma_k.  \]
\item[(d)] \[ r_{\delta - \la} = r_\la t_{-\gamma_1 - \gamma_2}. \]
\item[(e)] \[ r_{\delta-\la} (n' \La_0 + \rho) = n' \La_0 + \rho + \la - \delta = (-\La_0 + \la - \delta) + \wrh. \]
\end{enumerate}
\end{lemma}
Using Lemma \ref{Lem 2.4}, we can rewrite \eqref{2.20} as follows:
\[ \begin{aligned}
\hat{R} \, \ch L(-\La_0 + \bar{ \La}_2) & = -\frac{1}{q} \sum_{w \in W} \epsilon (w) w \sum_{\substack{\gamma \in Q^\vee \\ (\gamma|\bar{ \La}_1) < 0 \\ (\gamma | \bar{ \La}_{n'}) \in 2 \ZZ  }} t_\gamma r_{\bar{ \La}_2} t_{-\gamma_1-\gamma_2} (e^{-\La_0 + \bar{ \La}_2 + \wrh -\delta}). \\
& = \sum_{w \in W} \epsilon (w) w \sum_{\substack{ j_1, \ldots, j_{n'} \in \ZZ  \\ j_1 < 0 \\ \sum_k j_k \in 2 \ZZ }} t_{(-j_2-1)\gamma_1 + (-j_1-1)\gamma_2 + \sum_{k=3}^{n'} j_k \gamma_k } (e^{-\La_0 + \bar{ \La}_2 + \wrh}).
\end{aligned} \]
Replacing in the last expression $ -j_2 -1  $ by $ j_1 $ and $ -j_1 -1 $ by $ j_2, $ we obtain:
\[ \hat{R} \, \ch L(-\La_0 + \bar{ \La}_2) = \sum_{w \in W} \epsilon (w) w \sum_{\substack{ j_1, \ldots, j_{n'} \in \ZZ  \\ j_2 \geq 0 \\ \sum_k j_k \in 2 \ZZ }} t_{\sum_k j_k \gamma_k} (e^{-\La_0 + \bar{ \La}_2 + \wrh}). \]
Now, by \eqref{2.18}, claim (b) of Theorem \ref{Th2.2} follows.

\section{A character formula for the Deligne series modules}

In this section we prove the following simple theorem.

\begin{theorem}
\label{Th 3.1}
Let $ \fg $ be a simply laced Lie algebra of rank $ \ell $ (so that $ \al^\vee_i = \al_i $), and $ \La $ be a weight of $ \wg $ of level $ k \in \ZZ_{<0}, $ such that the following conditions hold:
\begin{enumerate}
\item[(i)] $ (\La | \al_i) \in \ZZ_{\geq 0} $ for $ i = 1, \ldots, \ell, $
\item[(ii)] there exists a positive root $ \al  $ of $ \fg,  $ such that $ (\La + \wrh | \delta - \al ) = 0, $
\item[(iii)] if $ \beta \in \hat{\Delta}_+ $ is orthogonal to $ \La + \wrh,  $ then $ \beta = \delta - \al, $
\item[(iv)] (extra hypothesis) in the character formula \eqref{0.09} one has:
\[ c(\gamma):= c(t_\gamma) = (\text{linear function in } \gamma \in Q) + \text{ const}. \]
\end{enumerate}
Then
\begin{equation}\label{3.01}
\hat{R} \, \ch \LLa = \half \sum_{w \in W} \epsilon (w) w \sum_{\gamma \in Q} ((\al|\gamma)+1) t_\gamma (e^{\La + \wrh}).
\end{equation}
\end{theorem}

Note that Theorem \ref{Th1.1} shows that the extra hypothesis fails for $ \fg = s\ell_n, \ k = -1. $ However, the comparison with \cite{Kaw15}, \cite{AK16} indicates that the following conjecture may hold.
\begin{conjecture}
\label{Conj 3.1}
If $ \fg = D_4, E_6, E_7,  $ or $ E_8, $ then the extra hypothesis (iv) holds (hence the character formula (\ref{3.01}) holds).
\end{conjecture}
\begin{remark}\label{rem 3.1}
  If $ \La = k \La_0 $ for $ k \in \ZZ_{<0}, $ conditions (i)--(iii) of Theorem \ref{Th 3.1} hold for $ k = ( - \frac{h^\vee}{6} -1 ) + s, $ where $ s = 0, 1, \ldots, b-1 $ and $ b = 2,3,4,6 $ for $ \fg = D_4, E_6, E_7, E_8 $ respectively.
  This explains the name ``Deligne series modules'', cf. \cite{Kaw15}, \cite{AK16}, \cite{AM16}.
\end{remark}


The proof of Theorem \ref{Th 3.1} is easy (but it is probably quite hard to verify the exta hypothesis (iv)).
Indeed, by (i) $ \ch \LLa $ is $ W $-invariant, and, by (ii), $k+h^\vee >0$, hence
\eqref{0.09} holds and it can be rewritten as 
\begin{equation}\label{3.02}
\hat{R} \, \ch \LLa = \sum_{w \in W} \epsilon (w) w \sum_{\gamma \in Q} c (\gamma) t_\gamma (e^{\La + \wrh}).
\end{equation}
Since, by (ii) we have $ r_{\delta - \al} (\La + \wrh) = \La + \wrh,  $ and also $ r_{\delta-\al}  = r_\al t_{-\al}, r_\al t_\gamma r_\al = t_{r_\al (\gamma)}, \epsilon (wr_\al) = - \epsilon (w), $, we can rewrite \eqref{3.02} as 
\[ \hat{R} \, \ch \LLa = - \sum_{w \in W} \epsilon (w) w \sum_{\gamma \in Q} c (\gamma) t_{r_\al (\gamma)-\al}  (e^{\La + \wrh}). \]
Replacing in this formula $ \gamma $ by $ r_\al (\gamma )- \al, $ we obtain
\begin{equation}\label{3.03}
\hat{R} \, \ch \LLa = - \sum_{w \in W} \epsilon (w) w \sum_{\gamma \in Q} c (r_\al (\gamma)-\al) t_\gamma (e^{\La + \wrh}).
\end{equation}
Taking half the sum of \eqref{3.02} and \eqref{3.03}, we obtain
\begin{equation}\label{3.04}
\hat{R} \, \ch \LLa = \half \sum_{w \in W} \epsilon (w) w \sum_{\gamma \in Q} \tilde{c} (\gamma) t_\gamma (e^{\La + \wrh})
\end{equation}
where $ \tilde{c} (\gamma) = c(\gamma) - c(r_\al (\gamma)-\al). $ The function $ \tilde{c} (\gamma) $ has the following two properties:
\begin{equation}\label{3.05}
\tilde{c} (\gamma) = -\tilde{c} (r_\al (\gamma)-\al), \gamma \in Q,
\end{equation}
\begin{equation}\label{3.06}
\tilde{c} (\gamma) = (\text{linear function in } \gamma) + \text{ const, } \gamma \in Q.
\end{equation}
By \eqref{3.06}, which holds due to the condition (iv), we can write for some $ \beta \in \fh^*, a \in \CC: $
\begin{equation}\label{3.07}
\tilde{c} (\gamma) = (\gamma | \beta) + a, \ \gamma \in Q.
\end{equation}
Then we obtain for all $ \gamma \in Q $:
\begin{equation}\label{3.08}
\tilde{c} (r_\al (\gamma)-\al) = (\gamma| \beta) - (\al | \beta) (\gamma|\al) - (\al | \beta) +a
\end{equation}
Since, by \eqref{3.05}, \eqref{3.07} = -\eqref{3.08}, we obtain:
\[ 2 (\gamma|\beta) - (\al|\beta) (\al|\al) - (\al|\beta) + 2a = 0, \text{ for all } \gamma \in Q. \]
Hence $ \beta = \half (\al|\beta)\al $ and $ a = \half (\al|\beta). $
Therefore, by \eqref{3.07}, we obtain
\begin{equation}\label{3.09}
\tilde{c} (\gamma)  = \text{const. } \times ((\gamma|\al)+1) \text{ for all } \gamma \in Q.
\end{equation}
By \eqref{3.04} and \eqref{3.09}, we have:
\begin{equation}\label{3.10}
\hat{R} \, \ch \LLa = \text{const. } \times \sum_{w \in W} \epsilon(w) w \sum_{\gamma \in Q} ((\gamma|\al)+1) t_\gamma (e^{\La + \wrh}).
\end{equation}
Since the stabilizer in $ \hat{W} $ of any $ \la \in \fh^* $ of positive level
is generated by reflections $ r_\al $, $\al \in \hat{\Delta}^{\text{re}}_+$ fixing $ \la $ \cite{K90},
by the conditions (ii) and (iii) we see that $ \hat{W}_{\La+\wrh} = \{ 1, r_{\delta-\al}\} $. It follows that const. $ = \half $ in \eqref{3.10}, proving \eqref{3.01}.

\begin{conjecture}
\label{Conj 3.4} 
If $ \fg = D_4, E_6, E_7 $ or $ E_8 $ and $ k = -1, -2, \ldots, -b, $ where $ b = 2,3,4, $ or 6 respectively, then all irreducible modules from the category $ \mathcal{O} $ of the vertex algebra $ L(k\La_0) $ are those from Theorem \ref{Th 3.1}. (It follows from \cite{AM16} that all these vertex algebras are quasilisse, hence, by \cite{AK16}, have only finitely many irreducible modules in the category $ \mathcal{O}. $)
\end{conjecture}
\begin{example}
\label{Ex 3.5}
Let $ \fg = D_4, \ k = -1. $ Then the following $ \La $'s satisfy the conditions (i), (ii), (iii) of Theorem \ref{Th 3.1} (we label the branching node of the Dynkin diagram of $D_4$ by 2):
\[ 
\begin{aligned}
& -\La_0; \ -2\La_0 + \La_i \ (i = 1,3,4); \ -3 \La_0 + \La_2; \\
& -3 \La_0 + \La_i + \La_j \ ((i,j) = (1,3), (1,4), (3,4)).
\end{aligned} \]
\end{example}
\begin{remark}
  \label{Rem 3.6}
Of course, one has a formula for homogeneous $q$-dimension, similar to that in Remark \ref{rem1.2} in all cases, considered in Sectins 2 and 3. We checked on the computer that in the case of $D_4$ it is compatible with the formula for $q$-dimension of $L(-2\La_0)$ from \cite{AK16}.
  \end{remark}
\begin{remark}  \label{Rem 3.7}
  Note that  $ ((\gamma|\alpha)+1)  t_\gamma (e^{\Lambda+\hat{\rho}})
  =\frac{1}{k+h^\vee}D_\alpha t_\gamma (e^{\Lambda+\hat{\rho}})$, where $D_\al$ is the derivative in the direction $\al$. Hence the RHS of \eqref{3.10} is a linear combination of derivatives of
  theta functions.
 \end{remark}
\begin{remark}
  \label{Rem 3.8}
  By Theorem 3.1 from \cite{KRW03} one has:
  \[\hat{R}(\fg,f) ch_{H(\Lambda)}(\tau, h)=(\hat{R}_nch_\Lambda)(\tau, -\tau x+h, \tau/4),
  \,\,\hbox{where}\,\, h\in \fh^f,\]
  for any $W$-algebra $W^k(\fg,f)$, obtained by the quantum Hamiltonian reduction
  of the $\hat{\fg}$-module $L(\La)$ of level $k$.
  Here $\hat{R}(\fg,f)$ is the denominator of $W^k(\fg,f)$,
  $\hat{R}_n=q^{\dim \fg/24}\hat{R}$ is the normalised affine Weyl denominator, $ch_\Lambda =q^{m_\Lambda}ch L(\Lambda) $
  is the normalized character \cite{K90}. In particular, if $\fg =D_4, E_6,E_7, E_8$, $k=-b$,
  $f$ is the minimal nilpotent element of $\fg$, and $L(\Lambda)$ is a module of level $k$ over the corresponding simple vertex algebra,
  then the simple $W$-algebra $W_k(\fg,f)$ is 1-dimensional \cite{AKMPP16}, hence $ch_{H(\Lambda)}=1$, and we get a formula relating $\hat{R}(\fg,f)$ to $ch L(\Lambda)$ for $z=-\tau x+h$.
\end{remark}

\end{document}